\pdfoutput=1
\RequirePackage{ifpdf}
\ifpdf 
\documentclass[pdftex]{sigma}
\else
\documentclass{sigma}
\fi

\usepackage{enumerate,young}
\usepackage[hcentermath]{youngtab}
\newlength\yStones
\newlength\xStones
\newlength\xxStones
 \DeclareMathOperator\Specm{Specm}
 
\DeclareMathOperator\Span{span} \DeclareMathOperator\Ker{Ker}
\DeclareMathOperator\Supp{Supp}

\numberwithin{equation}{section}

\newtheorem{Theorem}{Theorem}[section]
\newtheorem{Corollary}[Theorem]{Corollary}
\newtheorem{Lemma}[Theorem]{Lemma}
\newtheorem{Proposition}[Theorem]{Proposition}
{ \theoremstyle{definition}
\newtheorem{Definition}[Theorem]{Definition}

\newtheorem{Example}[Theorem]{Example}
\newtheorem{Notation}[Theorem]{Notation}
\newtheorem{Remark}[Theorem]{Remark} }

\makeatletter
\def\Stones{\pst@object{Stones}}
\def\Stones@i#1{%
  \pst@killglue%
  \begingroup%
  \use@par%
  \setlength\xxStones{\xStones}%
  \expandafter\Stones@ii#1,,\@nil
  \endgroup
  \mathfrak{gl}obal\addtolength\xStones{0.6cm}%
  \mathfrak{gl}obal\addtolength\yStones{-7.5mm}}%
\def\Stones@ii#1,#2,#3\@nil{%
  \rput(\xxStones,\yStones){%
    \psframebox[framesep=0]{%
      \parbox[c][6mm][c]{11mm}{\makebox[11mm]{$#1$}}}}%
  \addtolength\xxStones{1.2cm}%
  \ifx\relax#2\relax\else\Stones@ii#2,#3\@nil\fi}
\makeatother

\def\Stone#1{\fbox{\makebox[10mm]{\strut#1}}\kern2pt}

\begin{document}

\newcommand{\arXivNumber}{1409.8413}

\allowdisplaybreaks

\renewcommand{\thefootnote}{$\star$}

\renewcommand{\PaperNumber}{018}

\FirstPageHeading

\ShortArticleName{Irreducible Generic Gelfand--Tsetlin Modules of $\mathfrak{gl}(n)$}

\ArticleName{Irreducible Generic Gelfand--Tsetlin Modules of $\boldsymbol{\mathfrak{gl}(n)}$\footnote{This paper is
a~contribution to the Special Issue on New Directions in Lie Theory.
The full collection is available at
\href{http://www.emis.de/journals/SIGMA/LieTheory2014.html}{http://www.emis.de/journals/SIGMA/LieTheory2014.html}}}

\Author{Vyacheslav FUTORNY~$^\dag$, Dimitar GRANTCHAROV~$^\ddag$ and Luis Enrique RAMIREZ~$^\dag$}

\AuthorNameForHeading{V.~Futorny, D.~Grantcharov and L.E.~Ramirez}

\Address{$^\dag$~Instituto de Matem\'atica e Estat\'istica, Universidade de S\~ao Paulo, S\~ao Paulo SP, Brasil}
\EmailD{\href{mailto:futorny@ime.usp.br}{futorny@ime.usp.br}, \href{mailto:luiser@ime.usp.br}{luiser@ime.usp.br}}

\Address{$^\ddag$~Department of Mathematics, University of Texas at Arlington, Arlington, TX 76019, USA}
\EmailD{\href{grandim@uta.edu}{grandim@uta.edu}}

\ArticleDates{Received October 01, 2014, in f\/inal form February 24, 2015; Published online February 28, 2015}

\Abstract{We provide a~classif\/ication and explicit bases of tableaux of all irreducible generic Gelfand--Tsetlin modules
for the Lie algebra $\mathfrak{gl}(n)$.}

\Keywords{Gelfand--Tsetlin modules; Gelfand--Tsetlin basis; tableaux realization}

\Classification{17B67}

\renewcommand{\thefootnote}{\arabic{footnote}}
\setcounter{footnote}{0}

\section{Introduction}
Let $\mathfrak g$ be a~complex f\/inite-dimensional semisimple Lie algebra.
The category of weight modules of $\mathfrak g$ is interesting on its own on the one hand, and it contains some
fundamental subcategories like the category $\mathcal O$, categories of parabolically induced modules, Harish-Chandra
modules on the other.
A~weight $\mathfrak g$-module is a~module which is a~direct sum of simple $\mathfrak h$-modules, where~$\mathfrak h$ is
a~f\/ixed Cartan subalgebra of $\mathfrak g$.
The classif\/ication of the simple weight modules is a~very hard problem which is solved only for $\mathfrak g =
\mathfrak{sl}(2)$.
However, the classif\/ication of the simple objects is known for various subcategories of weight modules, including those
with f\/inite weight multiplicities~\cite{Fe, M}.

The classif\/ication of the simple weight $\mathfrak{sl}(2)$-modules involves two parameters that correspond to
eigenvalues of the generators of a~maximal commutative subalgebra of $U(\mathfrak{sl}(2))$, the {\it Gelfand--Tsetlin
subalgebra}.
Such subalgebra can be def\/ined for any $\mathfrak{sl}(n)$ and has a~joint spectrum on every f\/inite-dimensional module.
This observation leads naturally to the def\/inition of a~{\it Gelfand--Tsetlin module}: a~module that is the direct sum
of its common generalized eigenspaces with respect to the Gelfand--Tsetlin subalgebra~$\Gamma$.
Such modules were introduced in~\cite{DFO3, DFO2,DFO1}.
Note that an irreducible Gelfand--Tsetlin modules does not need to be~$\Gamma$-diagonalizable~\cite{Fut1}.

Gelfand--Tsetlin subalgebras and modules appear in various contexts.
Such subalgebras were considered in~\cite{Vi} in connection with subalgebras of maximal Gelfand--Kirillov dimension in
the universal enveloping algebra of a~simple Lie algebra.
Furthermore, Gelfand--Tsetlin subalgebras are related to: general hypergeometric functions on the complex Lie group
${\rm GL}(n)$~\cite{Gr1,Gr2}; solutions of the Euler equation~\cite{Vi}; and problems in classical mechanics in
general~\cite{KW-1,KW-2}.

One natural question is to attempt the classif\/ication of all irreducible Gelfand--Tsetlin mo\-du\-les of $\mathfrak{sl}(n)$.
An explicit construction of all irreducible Gelfand--Tsetlin modules for the case $n=3$ was recently obtained
in~\cite{FGR}.
Various partial results for $\mathfrak{sl}(3)$ were previously obtained in~\cite{Britten-Futorny-Lemire, Fut1,Fut3,Fut2,
Fut4}.

A {\it generic Gelfand--Tsetlin module} is a~module spanned by tableaux with noninteger dif\/fe\-ren\-ces of entries in each
row (see Def\/inition~\ref{def-gen}).
The present paper provides a~classif\/ication of all irreducible generic Gelfand--Tsetlin modules of $\mathfrak{sl}(n)$
extending the result in~\cite{Ramirez-theses} for $n=3$.
For simplicity we work with $\mathfrak{gl}(n)$ instead of~$\mathfrak{sl}(n)$.
We also obtain an explicit construction of all irreducible generic modules providing a~Gelfand--Tsetlin type basis.

The organization of the paper is as follows.
In Section~\ref{section3} we introduce some basic def\/initions and preparatory results on Gelfand--Tsetlin modules.
In Section~\ref{section4} we list the Gelfand--Tsetlin formulas and use them to recall the classical result of Gelfand
and Tsetlin for f\/inite-dimensional $\mathfrak{gl}(n)$-modules.
In Section~\ref{section5} we introduce the notion of generic Gelfand--Tsetlin module and recall the classif\/ication of
irreducible generic Gelfand--Tsetlin modules of~$\mathfrak{gl} (3)$.
The main theorem in the paper, the classif\/ication of irreducible generic Gelfand--Tsetlin $\mathfrak{gl} (n)$-modules,
is included in Section~\ref{section6}.
In the last section we compute the number of irreducible Gelfand--Tsetlin modules in the so-called generic blocks.

\section{Notation and conventions}

Throughout the paper we f\/ix an integer $n\geq 2$.
The ground f\/ield will be ${\mathbb C}$.
For $a \in {\mathbb Z}$, we write $\mathbb Z_{\geq a}$ for the set of all integers~$m$ such that $m \geq a$.
Similarly, we def\/ine $\mathbb Z_{< a}$, etc.
By~$\mathfrak{gl}(n)$ we denote the general linear Lie algebra consisting
of all $n\times n$ complex matrices, and by $\{E_{i,j}\,|\, 1\leq i,j \leq n\}$,
the standard basis of $\mathfrak{gl}(n)$ of elementary matrices.
We f\/ix the standard Cartan subalgebra $\mathfrak h$, the standard triangular decomposition and the corresponding basis
of simple roots of $\mathfrak{gl}(n)$.
The weights of $\mathfrak{gl}(n)$ will be written as~$n$-tuples $(\lambda_1,\dots,\lambda_n)$.

For a~Lie algebra ${\mathfrak a}$ by $U(\mathfrak a)$ we denote the universal enveloping algebra of ${\mathfrak a}$.
Throughout the paper $U = U(\mathfrak{gl} (n))$.
For a~commutative ring~$R$, by $\Specm R$ we denote the set of maximal ideals of~$R$.

We will write the vectors in $\mathbb{C}^{\frac{n(n+1)}{2}}$ in the following form:
\begin{gather*}
L=(l_{ij})=(l_{n1},\dots,l_{nn}| \dots|l_{21}, l_{22}|l_{11}).
\end{gather*}
For $1\leq j \leq i \leq n$, $\delta^{ij} \in {\mathbb Z}^{\frac{n(n+1)}{2}}$ is def\/ined by $(\delta^{ij})_{ij}=1$ and
all other $(\delta^{ij})_{k\ell}$ are zero.

For $i>0$ by $S_i$ we denote the~$i$th symmetric group.
Throughout the paper we set $G:=S_n\times\dots \times S_1$.

\section{Gelfand--Tsetlin modules}
\label{section3}

Recall that $U=U(\mathfrak{gl} (n))$.
Let for $m\leqslant n$, $\mathfrak{gl}_{m}$ be the Lie subalgebra of $\mathfrak{gl} (n)$ spanned by $\{E_{ij} \,|\,
i,j=1,\ldots,m \}$.
We have the following chain
\begin{gather*}
\mathfrak{gl}_1\subset \mathfrak{gl}_2\subset \cdots \subset \mathfrak{gl}_n.
\end{gather*}
It induces the chain $U_1\subset$ $U_2\subset \dots \subset U_n$ for the universal enveloping algebras
$U_{m}=U(\mathfrak{gl}_{m})$, $1\leq m\leq n$.
Let $Z_{m}$ be the center of $U_{m}$.
The subalgebra of~$U$ generated by $\{Z_m \,|\, m=1,\ldots, n \}$ will be called the \emph{$($standard$)$ Gelfand--Tsetlin
subalgebra} of~$U$ and will be denoted by ${{\Gamma}}$~\cite{DFO3}.

\begin{Definition}
A~f\/initely generated~$U$-module~$M$ is called a~\emph{Gelfand--Tsetlin module $($with respect to ${\Gamma})$} if
\begin{gather*}
M=\bigoplus_{\mathsf m\in\Specm{\Gamma}}M(\mathsf m),
\end{gather*}
where $M(\mathsf m)=\{v\in M\,|\, \mathsf m^{k}v=0~\text{for some}~k\geq 0\}$.
\end{Definition}

For each $\mathbf m\in\Specm{\Gamma}$ we have associated a~character $\chi_{\mathbf m}:{\Gamma}\rightarrow{\Gamma}/\mathbf
m\sim\mathbb{C}$.
In the same way, for each non-zero character $\chi:{\Gamma}\rightarrow\mathbb{C}$ we have that $Ker(\chi)$ is a~maximal
ideal of ${\Gamma}$.
So, we have a~natural identif\/ication between characters of ${\Gamma}$ and elements of $\Specm{\Gamma}$.
Using characters we can def\/ine Gelfand--Tsetlin modules.
A~$U$-module~$M$ is called \emph{Gelfand--Tsetlin module} (with respect to~$\Gamma$) if
\begin{gather*}
M=\bigoplus_{\chi\in\Gamma^{*}}M(\chi),
\end{gather*}
where $M(\chi)=\{v\in M: \forall\, g\in\Gamma, \, \exists\, k\in\mathbb{Z}_{>0}~\text{such that}~(g-\chi(g))^{k}v=0 \}$.
The {\it Gelfand--Tsetlin support} of~$M$ is the set $\Supp_{\rm GT}(M):=\{\chi\in{\Gamma}^{*}:M(\chi)\neq 0\}$.

\begin{Lemma}
\label{lemma elements of Gamma separates tableaux}
Any submodule of a~Gelfand--Tsetlin module over $\mathfrak{gl}(n)$ is a~Gelfand--Tsetlin mo\-dule.
\end{Lemma}

\begin{proof}
The proof is standard, but for a~sake of completeness, we provide the important details.
Let~$M$ be a~Gelfand--Tsetlin $\mathfrak{gl}(n)$-module and~$N$ any submodule of~$M$.
We will prove that, if $\{\chi_{1},\dots,\chi_{k} \}$ is a~set of distinct Gelfand--Tsetlin characters in
$\Supp_{\rm GT}(M)$ such that $\sum\limits_{i=1}^{k}v_{i}\in N$ with $v_{i}\in M(\chi_{i})$, then $v_{i}\in N$ for all
$i=1,\ldots,k$.

Without loss of generality we assume that $k=2$.
Since $\chi_{1}\neq\chi_{2}$, there exist $g\in\Gamma$ and $r\leq s$ in $\mathbb{Z}_{\geq 0}$ such that
$\chi_{1}(g)\neq\chi_{2}(g)$, $(g-\chi_{1}(g))^{r}(v_{1})=0$ and $(g-\chi_{2}(g))^{s}(v_{2})=0$.
Let $a:=\chi_{1}(g)$ and $b:= \chi_{2}(g)$, Then, if $w=v_{1}+v_{2}$ we have $(g-b)^{s}w=(g-b)^{s}v_{1}\in N$.
Let $y:= (g-b)^sv_1$.
We have that $y \in N$ on one hand and
\begin{gather*}
y=((g-a)+(a-b))^{s}v_{1}=\sum\limits_{k=0}^{r-1} {s\choose k}(a-b)^{s-k}(g-a)^{k}v_{1}\in N
\end{gather*}
on the other.
As ${s\choose k}(a-b)^{s-k}\neq 0$ for any~$k$, using that $(g-a)^{r-1}y\in N$, we obtain \mbox{$(g-a)^{r-1}v_1\in N$}.
Reasoning in the same way, from $(g-a)^{r-i}y\in N$, and $(g-a)^{r-1}v_{1},\ldots$, $(g-a)^{r-i+1}v_{1}\in N $ we obtain $
x^{r-i}v_{1}\in N$.
Hence $v_{1}\in N$ and consequently, $v_{2}\in N$.
\end{proof}

One can choose the following generators of~$\Gamma$: $\{c_{mk} \,|\, 1\leq k\leq m\leq n \}$, where
\begin{gather}
\label{equ_3}
c_{mk}= {\sum\limits_{(i_1,\ldots,i_k)\in \{1,\ldots,m \}^k}} E_{i_1 i_2}E_{i_2 i_3}\cdots E_{i_k i_1}.
\end{gather}

Let~$\Lambda$ be the polynomial algebra in the variables $\{\lambda_{ij}\, | \, 1\leqslant j\leqslant i\leqslant n \}$.
The action of the symmetric group $S_i$ on $\{\lambda_{ij} |$ $1\leqslant j\leqslant i\}$ induces the action of $G = S_n
\times\dots \times S_1$ on~$\Lambda$.
There is a~natural embedding $\imath:{{\Gamma}}{\longrightarrow}$~$\Lambda$ given~by
$\imath(c_{mk})=\gamma_{mk}(\lambda)$, where
\begin{gather}
\label{def-gamma}
\gamma_{mk}(\lambda)=\sum\limits_{i=1}^m (\lambda_{mi}+m-1)^k \prod\limits_{j\ne i} \left(1-\frac{1}{\lambda_{mi}-\lambda_{mj}} \right).
\end{gather}
Hence,~$\Gamma$ can be identif\/ied with $G$-invariant polynomials in~$\Lambda$.

\begin{Remark}
\label{maximal ideals differ by permutations}
In what follows, we will identify the set $\Specm \Lambda$ of maximal ideals of~$\Lambda$ with the set
$\mathbb{C}^{\frac{n(n+1)}{2}}$.
Then we have a~surjective map $\pi: \Specm\Lambda \rightarrow \Specm{\Gamma}$.
Moreover, since~$\Lambda$ is integral over~$\Gamma$, there are f\/initely many maximal ideals of~$\Lambda$ that map to
a~f\/ixed maximal ideal of~$\Gamma$.
The dif\/ferent maximal ideals of~$\Lambda$ are obtained from each other under permutations in the group~$G$.
\end{Remark}
If $\pi(\ell)=\mathsf m$ for some $\ell\in \Specm \Lambda$, then we write $\ell=\ell_{\mathsf m}$ and say that
$\ell_{\mathsf m}$ is {\it lying over} $\mathsf m$.

\section[Finite-dimensional modules of $\mathfrak{gl}(n)$]{Finite-dimensional modules of $\boldsymbol{\mathfrak{gl}(n)}$}
\label{section4}

In this section we recall a~classical result of Gelfand and Tsetlin which provides an explicit basis for every
irreducible f\/inite-dimensional $\mathfrak{gl}(n)$-module.

\begin{Definition}
For a~vector $L=(l_{ij})$ in $\mathbb{C}^{\frac{n(n+1)}{2}}$, by $T(L)$ we will denote the following array with entries
$\{l_{ij}:1\leq j\leq i\leq n\}$

\begin{center}
\Stone{\mbox{\scriptsize {$l_{n1}$}}}\Stone{\mbox{\scriptsize {$l_{n2}$}}}\hspace{1cm} $\dots$ \hspace{1cm}
\Stone{\mbox{\scriptsize {$l_{n,n-1}$}}}\Stone{\mbox{\scriptsize {$l_{nn}$}}}
\\[0.2pt] \Stone{\mbox{\scriptsize {$l_{n-1,1}$}}}\hspace{1.5cm} $\dots$ \hspace{1.5cm} \Stone{\mbox{\tiny
{$l_{n-1,n-1}$}}}
\\[0.3cm] \hspace{0.2cm}$\dots$ \hspace{0.8cm} $\dots$ \hspace{0.8cm} $\dots$
\\[0.3cm] \Stone{\mbox{\scriptsize {$l_{21}$}}}\Stone{\mbox{\scriptsize {$l_{22}$}}}
\\[0.2pt] \Stone{\mbox{\scriptsize {$l_{11}$}}}
\end{center}
Such an array will be called a~\emph{Gelfand--Tsetlin tableau} of height~$n$.
A~Gelfand--Tsetlin tableau of height~$n$ is called \emph{standard} if $l_{ki}-l_{k-1,i}\in\mathbb{Z}_{\geq 0}$ and
$l_{k-1,i}-l_{k,i+1}\in\mathbb{Z}_{> 0}$ for all $1\leq i\leq k\leq n-1$.
\end{Definition}
Note that, for sake of convenience, the second condition above is slightly dif\/ferent from the original condition
in~\cite{GT}.

\begin{Theorem}[\cite{GT}]\label{Gelfand--Tsetlin theorem}
Let $L(\lambda)$ be the finite-dimensional irreducible module over $\mathfrak{gl}(n)$ of highest weight
$\lambda=(\lambda_{1},\ldots,\lambda_{n})$.
Then there exist a~basis of $L(\lambda)$ consisting of all standard tableaux $T(L) = T(l_{ij})$ with fixed top row
$l_{nj}=\lambda_j-j+1$.
Moreover, the action of the generators of~$\mathfrak{gl}(n)$ on~$L(\lambda)$ is given by the Gelfand--Tsetlin
formulas:
\begin{gather}
E_{k,k+1}(T(L)) =-\sum\limits_{i=1}^{k}\left(\frac{\prod\limits_{j=1}^{k+1}(l_{ki}-l_{k+1,j})}{\prod\limits_{j\neq
i}^{k}(l_{ki}-l_{kj})}\right)T\big(L+\delta^{ki}\big),
\nonumber
\\
E_{k+1,k}(T(L)) =\sum\limits_{i=1}^{k}\left(\frac{\prod\limits_{j=1}^{k-1}(l_{ki}-l_{k-1,j})}{\prod\limits_{j\neq
i}^{k}(l_{ki}-l_{kj})}\right)T\big(L-\delta^{ki}\big),
\nonumber
\\
E_{kk}(T(L)) =\left(k-1+\sum\limits_{i=1}^{k}l_{ki}-\sum\limits_{i=1}^{k-1}l_{k-1,i}\right)T(L),
\label{Gelfand--Tsetlin formulas}
\end{gather}
if the new tableau $T(L\pm\delta^{ki})$ is not standard, then the corresponding summand of $E_{k,k+1}(T(L))$ or
$E_{k+1,k}(T(L))$ is zero by definition.
Furthermore, for $s\leq r$,
\begin{gather}
\label{action of Gamma in finite dimensional modules}
c_{rs}(T(L))=\gamma_{rs}(l)T(L),
\end{gather}
where $\{c_{rs}\}$ are the generators of~$\Gamma$ defined in~\eqref{equ_3} and $\gamma_{rs}$ are defined
in~\eqref{def-gamma} $($see~{\rm \cite{Zh})}.
\end{Theorem}
\noindent The formulas above are called \emph{Gelfand--Tsetlin formulas} for $\mathfrak{gl}(n)$.
These formulas were extended to the case of $U_{q}(\mathfrak{gl}(n))$ in~\cite{Tar-Naz}.

\section[Generic Gelfand--Tsetlin modules of $\mathfrak{gl}(n)$]{Generic Gelfand--Tsetlin modules of $\boldsymbol{\mathfrak{gl}(n)}$}
\label{section5}

Theorem~\ref{Gelfand--Tsetlin theorem} gives an explicit realization of any irreducible f\/inite-dimensional
$\mathfrak{gl}(n)$-module.
Using the Gelfand--Tsetlin formulas, Drozd, Futorny and Ovsienko def\/ined the class of inf\/inite-dimensional generic
modules for $\mathfrak{gl}(n)$ in~\cite{DFO3}.

\begin{Definition}
\label{def-gen}
A~Gelfand--Tsetlin tableau $T(L)$ \big(equivalently, $L\in{\mathbb C}^{\frac{n(n+1)}{2}}$\big) is called \emph{generic} if
$l_{ki}-l_{kj}\notin\mathbb{Z}$ for all $1\leq i\neq j \leq k \leq n-1$.
A~character~$\chi$ and $\mathsf n = \Ker \chi$ are called \emph{generic} if $\ell_{\mathsf n}$ is generic for one choice
(hence for all choices) of $\ell_{\mathsf n}$ lying over $\mathsf n$.
A~Gelfand--Tsetlin module~$M$ will be called a~\emph{generic Gelfand--Tsetlin module} if every $\mathsf n$ in
$\Supp_{\rm GT}(M)$ is generic.
\end{Definition}

\begin{Theorem}[\protect{\cite[Section 2.3]{DFO3} and~\cite[Theorem 2]{Maz2}}]
Let $T(L) = T(l_{ij})$ be a~generic Gelfand--Tsetlin tableau of height~$n$.
Denote by ${\mathcal B}(T(L))$ the set of all Gelfand--Tsetlin tableaux $T(R) = T(r_{ij})$ satisfying $r_{nj}=l_{nj}$,
$r_{ij}-l_{ij}\in\mathbb{Z}$ for
$1\leq j\leq i \leq n-1$.
\begin{itemize}\itemsep=0pt
\item[$(i)$] The vector space $V(T(L)) = \Span {\mathcal B}(T(L))$ has a~structure of a~$\mathfrak{gl}(n)$-module with
action of the generators of $\mathfrak{gl}(n)$ given by the Gelfand--Tsetlin formulas~\eqref{Gelfand--Tsetlin formulas}.
\item[$(ii)$] The action of the generators of~$\Gamma$ on the basis elements of $V(T(L))$ is given by~\eqref{action of Gamma in finite dimensional modules}.
\item[$(iii)$] The $\mathfrak{gl}(n)$-module $V(T(L))$ is a~Gelfand--Tsetlin module all of whose Gelfand--Tsetlin
multiplicities are~$1$.
\end{itemize}
\end{Theorem}

\begin{Remark}

The basis of the module in the previous theorem is
\begin{gather*}
{\mathcal B}(T(L)) = \big\{T(L+z):z\in\mathbb{Z}^{\frac{n(n+1)}{2}}~\text{and}~z_{n1}=\cdots=z_{nn}=0\big\}.
\end{gather*}
By a~slight abuse of notation we will identify elements in $\mathbb{Z}^{\frac{n(n-1)}{2}}$ with elements
$z\in\mathbb{Z}^{\frac{n(n+1)}{2}}$ such that $z_{n1}=\cdots=z_{nn}=0$.
This will allow us to write $T(L+z)$ for $z\in\mathbb{Z}^{\frac{n(n-1)}{2}}$.
\end{Remark}

\begin{Remark}
\label{elements of Gamma separates tableaux}
In what follows, we will apply Lemma~\ref{lemma elements of Gamma separates tableaux} and use that the elements
of~$\Gamma$ separate the tableaux in the submodules of $V(T(L))$ in the following sense.
Let~$N$ be a~$\mathfrak{gl}(n)$-submodule of $V(T(L))$, $g \in \mathfrak{gl}(n)$, and $T(R)$ be a~tableau in~$N$.
Then, if $g \cdot T(R) = \sum\limits_i c_i T(R_i)$ for some distinct tableaux $T(R_i)$ in ${\mathcal B}(T(L))$ and
nonzero $c_i \in {\mathbb C}$, we have $T(R_i) \in N$ for all~$i$.
\end{Remark}

\begin{Theorem}
\label{uniqueness in generic case}
If $\mathsf n\in\Specm\Gamma$ is generic, then there exists a~unique irreducible Gelfand--Tsetlin module~$N$ such that
$N(\mathsf n)\neq 0$.
\end{Theorem}
\begin{proof}
Let $X_{\mathsf n}=U/U \mathsf n$.
We know that $X_{\mathsf n}=U/U \mathsf n$ is a~Gelfand--Tsetlin module.
Furthermore, any irreducible Gelfand--Tsetlin module~$M$ with $M(\mathsf n)\neq 0$ is a~homomorphic image of $X_{\mathsf
n}$, and~$X_{\mathsf n}(\mathsf n)$ maps onto~$M(\mathsf n)$.
Since both spaces $X_{\mathsf n}(\mathsf n)$ and $M(\mathsf n)$ are~$\Gamma$-modules then the projection $X_{\mathsf
n}(\mathsf n) \to M(\mathsf n)$ is a~homomorphism of~$\Gamma$-modules (see also~\cite[Corollary~5.3]{FO2}).
Taking into account that $\dim X_{\mathsf n}(\mathsf n) \leq 1$, we conclude that $X_{\mathsf n}$ has a~unique maximal
submodule (which does not intersect $X_{\mathsf n}(\mathsf n)$) and hence there exist a~unique irreducible module~$N$
with $N(\mathsf n)\neq 0$.
\end{proof}

\begin{Definition}
\label{def-cont}
If $T(R)$ is a~generic tableau and $\mathsf r\in\Specm\Gamma$ corresponds to~$R$ then, the unique module~$N$ such that
$N(\mathsf r) \neq 0$ is called the \emph{irreducible Gelfand--Tsetlin module containing $T(R)$}, or simply, the
\emph{irreducible module containing $T(R)$}.
\end{Definition}

Our goal is to describe explicitly the irreducible Gelfand--Tsetlin module containing $T(R)$ for every generic tableau $T(R)$.
Below we recall how this is achieved in the case $n=3$ in~\cite{Ramirez}.
One should note that the methods used in~\cite{Ramirez} involve direct computations based on a~case-by-case
consideration, while in the present paper we provide an invariant proof.
Also, we reformulate the result in~\cite{Ramirez} in terms of $T(L+z)$.

For any tableau $T(R)\in\{T(L+z):z\in \mathbb{Z}^{3}\}$ and any $1<p\leq 3$, $1\leq s\leq p$, and $1\leq u\leq p-1$, def\/ine
\begin{gather*}
\Omega^{+}(T(R)) :=\{(p,s,u):r_{p,s}-r_{p-1,u}\in\mathbb{Z}_{\geq 0}\}.
\end{gather*}

\begin{Theorem}[\cite{Ramirez}]
\label{Basis for irreducible generic modules gl(3)}
If $T(L)$ is a~generic Gelfand--Tsetlin tableau of height~$3$, then the following is a~basis for the irreducible
$\mathfrak{gl}(3)$-module containing~$T(L)$:
\begin{gather*}
\mathcal{I}(T(L)):=\big\{T(L+z):z\in\mathbb{Z}^{3}~\text{and}~\Omega^{+}(T(L))=\Omega^{+}(T(L+z))\big\}.
\end{gather*}
The action of $\mathfrak{gl}(3)$ on this irreducible module is given by the Gelfand--Tsetlin formulas.
\end{Theorem}

\begin{Example}
Consider $a,b,c\in\mathbb{C}$ such that $\{a-b,a-c,b-c\}\bigcap\mathbb{Z}=\varnothing$, $L=(a,b,c|a,b+1|a)$ and
\begin{center}
\hspace{1.5cm} \Stone{a}\Stone{b}\Stone{c}
\\[0.2pt] $T(L)$=\hspace{0.5cm} \Stone{a}\Stone{b+1}\hspace{1.2cm}
\\[0.2pt] \hspace{1.4cm} \Stone{a}
\\
\end{center}
then $\Omega^{+}(T(L))=\{(3,1,1), (2,1,1)\}$.
So, by Theorem~\ref{Basis for irreducible generic modules gl(3)}, the irreducible module con\-tai\-ning~$T(L)$ has basis
\begin{gather*}
\mathcal{I}(T(L))=\big\{T(L+(m,n,k)):(m,n,k)\in\mathbb{Z}^{3}, \ m\leq 0,\ k\leq m,~\text{and}~n>-1\big\}.
\end{gather*}
\end{Example}

\section[Classif\/ication of irreducible generic Gelfand--Tsetlin $\mathfrak{gl}(n)$-modules]{Classif\/ication of irreducible generic\\
Gelfand--Tsetlin $\boldsymbol{\mathfrak{gl}(n)}$-modules}
\label{section6}

In this section we prove the main result in the paper, i.e.~the generalization of Theorem~\ref{Basis for irreducible
generic modules gl(3)} for~$\mathfrak{gl}(n)$.
For convenience we introduce and recall some notation.

\begin{Notation}
\label{notation}
Let $T(L)=T(l_{ij})$ be a~f\/ixed tableau of height~$n$.
\begin{itemize}\itemsep=0pt
\item[(i)] $\mathcal{B}(T(L)):=\big\{T(L+z): z\in \mathbb{Z}^{\frac{n(n-1)}{2}}\big\}$.
\item[(ii)] $V(T(L)):= \Span \mathcal{B}(T(L))$.
\item[(iii)] For any $T(R)=T(r_{ij})\in \mathcal{B}(T(L))$ and for any $1<p\leq n$, $1\leq s\leq p$ and $1\leq u\leq
p-1$ we def\/ine:
\begin{itemize}\itemsep=0pt
\item[(a)] $\omega_{p,s,u}(T(R)):=r_{p,s}-r_{p-1,u}$;
\item[(b)] $\Omega(T(R)):=\{(p,s,u): \omega_{p,s,u}(T(R))\in \mathbb{Z}\}$;
\item[(c)] $\Omega^{+}(T(R)):=\{(p,s,u): \omega_{p,s,u}(T(R))\in \mathbb{Z}_{\geq 0}\}$;
\item[(d)] $\mathcal{N}(T(R)):=\{T(Q)\in \mathcal{B}(T(L)): \Omega^{+}(T(R))\subseteq \Omega^{+}(T(Q))\}$;
\item[(e)] $ W(T(R)):=\Span \mathcal{N}(T(R))$;
\item[(f)] $U\cdot T(R)$: the $\mathfrak{gl}(n)$-submodule of $V(T(L))$ generated by $T(R)$.
\end{itemize}
\end{itemize}
\end{Notation}

\subsection{Basis for the module generated by a~single tableau}
In order to f\/ind an explicit basis of every irreducible generic module, we f\/irst f\/ind a~basis of $U\cdot T(R)$ for any
tableau $T(R)$ in $\mathcal{B}(T(L))$.
\begin{Proposition}
\label{submodules of V(T(L))}
For any $T(R)\in\mathcal{B}(T(L))$, the Gelfand--Tsetlin formulas endow $ W(T(R))$ with a~$\mathfrak{gl}(n)$-module
structure.
\end{Proposition}
\begin{proof}
It is enough to prove $U\cdot T(Q)\subseteq W(T(R))$ for any $T(Q)=T(q_{ij})\in\mathcal{N}(T(R))$.
We will show $g\cdot T(Q)$ is in $W(T(R))$ for every (standard) generator~$g$ of $\mathfrak{gl}(n)$.

Suppose $g=E_{k,k+1}$ for some $1\leq k\leq n-1$.
By the Gelfand--Tsetlin formulas, we have
\begin{gather*}
E_{k,k+1}(T(Q))=-\sum\limits_{i=1}^{k}\left(\frac{\prod\limits_{j=1}^{k+1}(q_{ki}-q_{k+1,j})}{\prod\limits_{j\neq
i}^{k}(q_{ki}-q_{kj})}\right)T\big(Q+\delta^{ki}\big).
\end{gather*}
If $E_{k,k+1}(T(Q))\notin W(T(R))$, then there exist~$k$ and~$i$ such that $T(Q)\in \mathcal{N}(T(R))$ but
$T(Q+\delta^{ki})$ $\notin\mathcal{N}(T(R))$.
That implies
\begin{gather*}
\Omega^{+}(T(R))\subseteq\Omega^{+}(T(Q))~\text{and}~\Omega^{+}(T(R))\nsubseteq\Omega^{+}\big(T\big(Q+\delta^{ki}\big)\big).
\end{gather*}
Hence, there exists $(p,s,u)\in\Omega^{+}(T(R))$ such that $\omega_{p,s,u}(T(Q))\in\mathbb{Z}_{\geq 0}$ and
$\omega_{p,s,u}(T(Q+\delta^{ki}))$ $\notin\mathbb{Z}_{\geq 0}$.
The latter holds only in two cases:
\begin{gather*}
(p,s,u)\in \{(k,i,u), (k+1,s,i): 1\leq u\leq k-1; \, 1\leq s\leq k+1\}.
\end{gather*}
Note that if neither of these two cases hold, we have $\omega_{p,s,u}(T(Q+\delta^{ki}))=\omega_{p,s,u}(T(Q))$.
We consider now each of the two cases separately.
\begin{enumerate}[(i)]\itemsep=0pt
\item Suppose $(p,s,u)=(k,i,u)$.
Then $\omega_{k,i,u}(T(Q))=q_{ki}-q_{k-1,u}\in\mathbb{Z}_{\geq 0}$ and
$\omega_{k,i,u}(T(Q+\delta^{ki}))=(q_{ki}+1)-q_{k-1,u}\notin \mathbb{Z}_{\geq 0}$, which is impossible.
\item Suppose $(p,s,u)=(k+1,s,i)$.
Then
\begin{gather*}
\omega_{k+1,s,i}(T(Q))=q_{k+1,s}-q_{ki}\in\mathbb{Z}_{\geq 0}
\end{gather*}
{and}
\begin{gather*}
\omega_{k+1,s,i}(T(Q+\delta^{ki}))=q_{k+1,s}-(q_{ki}+1)\notin \mathbb{Z}_{\geq 0}.
\end{gather*}
Hence $q_{k+1,s}-q_{k,i}=0$ and then the coef\/f\/icient of $T(Q+\delta^{ki})$ in the decomposition of $E_{k,k+1}(T(Q))$ is
$-\frac{\prod\limits_{j=1}^{k+1}(q_{ki}-q_{k+1,j})}{\prod\limits_{j\neq i}^{k}(q_{ki}-q_{kj})}=0$.
\end{enumerate}
Therefore, the tableaux that appear with nonzero coef\/f\/icients in $E_{k,k+1}(T(Q))$ are elements of $N(T(R))$.
Hence, $E_{k,k+1}(T(Q))\in W(T(R))$.
The proof that $E_{k+1,k}(T(Q))\in W(T(R))$ is analogous to the one of $E_{k,k+1}(T(Q))$ $\in W(T(R))$.
The case $g=E_{kk}$ is trivial because $E_{kk}$ acts as a~multiplication by a~scalar on $T(Q)$ and
$T(Q)\in\mathcal{N}(T(R))\subseteq W(T(R))$.
\end{proof}

Given any tableau $T(R)$, there are three modules containing $T(R)$: $V(T(L))$, $ W(T(R))$ and $U\cdot T(R)$.
We will show that $ W(T(R))=U\cdot T(R)$.
For this we need the following lemmas.

\begin{Lemma}
\label{paths for connected tableaux}
Let $T(L)$ be a~generic tableau.
If $0\neq z\in\mathbb{Z}^{\frac{n(n-1)}{2}}$ is such that $\Omega^{+}(T(L))\subseteq\Omega^{+}(T(L+z))$ then, there
exist $i$, $j$ such that $z_{ij}\neq 0$ and
\begin{gather}
\label{inclutions of Omega +}
\Omega^{+}(T(L))\subseteq\Omega^{+}\big(T\big(L+z_{ij}\delta^{ij}\big)\big)\subseteq\Omega^{+}(T(L+z)).
\end{gather}
\end{Lemma}
\begin{proof}
We will use the following def\/inition in the proof of the lemma.

\begin{Definition}
Given a~generic tableau $T(R)\in\mathcal{B}(T(L))$, a~\emph{chain in $T(R)$ of length~$\ell$ starting in row~$d$} is
a~subset of the entries of $T(R)$, $C=\{r_{d-i,s^{(d-i)}}\}_{i=0,\ldots, \ell}$, where $1\leq s^{(d-i)}\leq d-i$ are
such that $r_{d-i,s^{(d-i)}}-r_{d-i-1,s^{(d-i-1)}}\in\mathbb{Z}$ for any $i=0,\ldots,\ell-1$
(i.e.~$\{(d-i,s^{(d-i)},s^{(d-i-1)})\}_{i=0,\ldots,\ell}$ $\subseteq\Omega(T(R))$). The chain is called maximal if
\begin{itemize}\itemsep=0pt
\item[(i)] $(d+1,i,s^{(d)})\notin\Omega(T(R))$ for any $1\leq i\leq d+1$,
\item[(ii)] $(d-\ell,s^{(d-\ell)},j)\notin\Omega(T(R))$ for any $1\leq j\leq d-\ell-1$.
\end{itemize}
\end{Definition}

For every $T(R)$ in $\mathcal{B}(T(L))$ we have that $\Omega^{+}(T(R))= \bigsqcup_{1\leq c\leq n}\Omega^{+}_{c}(T(R))$,
where $\Omega^{+}_{c}(T(R)):=\{(p,s,u)\in\Omega^{+}(T(R)): p=c\}$.
In particular,~\eqref{inclutions of Omega +} holds if and only if
\begin{gather}
\label{inclutions for Omega c +}
\Omega^{+}_{c}(T(L))\subseteq\Omega^{+}_{c}\big(T\big(L+z_{ij}\delta^{ij}\big)\big)\subseteq\Omega^{+}_{c}(T(L+z))
\end{gather}
for any $1\leq c\leq n$.
For $c\notin\{i,i+1 \}$ we have $\Omega^{+}_{c}(T(L))=\Omega^{+}_{c}(T(L+z_{ij}\delta^{ij}))$.
So, in order to verify~\eqref{inclutions for Omega c +}, it is enough to consider the cases $c=i, i+1$.

Lets consider $k$, $l$ such that $z_{kl}\neq 0$.
Set for convenience $Q:=L+z$.
There exists a~maximal chain~$C$ in $T(Q)$ of length~$\ell$, starting in row~$d$ such that $q_{kl}\in C$.
Suppose that $C=\{q_{[i]}\}_{i=0,\ldots,\ell}$ where $[i]:=(d-i,s^{(d-i)})$.
If $\ell=0$, then $C=\{q_{kl}\}$ and~\eqref{inclutions of Omega +} is obvious for $z_{ij}=z_{kl}$.

Let~$a$ and~$b$ be the minimum and maximum of $\{i:z_{[i]}\neq 0\}$, respectively.
We have
\begin{gather}
\Omega^{+}_{d-a+1}\big(T\big(L+z_{[a]}\delta^{[a]}\big)\big)=\Omega^{+}_{d-a+1}(T(L+z)),\nonumber
\\
\Omega^{+}_{d-b}\big(T\big(L+z_{[b]}\delta^{[b]}\big)\big)=\Omega^{+}_{d-b}(T(L+z)).\label{a and b behavior}
\end{gather}

Therefore~\eqref{inclutions for Omega c +} holds for the pairs $c=d-a+1$, $z_{ij}=z_{[a]}$ and $c=d-b$, $z_{ij}=z_{[b]}$,
respectively.
Now, let $a\leq m\leq b$ and consider the $4$ cases depending on what the signs of $z_{[a]}$ and $z_{[a+1]}$ are.
\begin{enumerate}[(i)]\itemsep=0pt
\item $z_{[m]}>0$ and $z_{[m+1]}\leq 0$.
In this case~\eqref{inclutions for Omega c +} holds for $c=d-m$ and $z_{ij}=z_{[m]}$.
In particular, if $z_{[a]}>0$ and $z_{[a+1]}\leq 0$, using the f\/irst equation in~\eqref{a and b behavior}, we conclude
that~\eqref{inclutions of Omega +} holds for $z_{ij}=z_{[a]}$.
\item $z_{[m]}<0$ and $z_{[m-1]}\geq 0$.
In this case~\eqref{inclutions for Omega c +} holds for $c=d-m+1$ and $z_{ij}=z_{[m-1]}$.
In particular, if $z_{[b]}<0$ and $z_{[b-1]}\geq 0$, using the second equation in~\eqref{a and b behavior} we conclude
that~\eqref{inclutions of Omega +} holds for $z_{ij}=z_{[b]}$.
\item $z_{[m]}>0$ and $z_{[m+1]}> 0$.
In this case~\eqref{inclutions for Omega c +} holds for $c=d-m$ and
\begin{gather*}
z_{ij}=
\begin{cases}
z_{[m]} & \text{if} \ \ l_{[m]}-l_{[m+1]}\in\mathbb{Z}_{\geq 0},
\\
z_{[m+1]} & \text{if} \ \ l_{[m+1]}-l_{[m]}\in\mathbb{Z}_{> 0}.
\end{cases}
\end{gather*}
\item $z_{[m]}<0$ and $z_{[m-1]}< 0$.
In this case~\eqref{inclutions for Omega c +} holds for $c=d-m+1$ and
\begin{gather*}
z_{ij}=
\begin{cases}
z_{[m]} & \text{if} \ \  l_{[m-1]}-l_{[m]}\in\mathbb{Z}_{\geq 0},
\\
z_{[m-1]} & \text{if} \ \ l_{[m]}-l_{[m-1]}\in\mathbb{Z}_{> 0}.
\end{cases}
\end{gather*}
\end{enumerate}

Now combining (i)--(iv) we reduce the proof to the following two cases:
\begin{enumerate}[(a)]\itemsep=0pt
\item $z_{[a]}>0, z_{[a+1]}>0, \ldots, z_{[b]}>0$ and for any $t=1,\ldots,b-a$,~\eqref{inclutions for Omega c +} holds
for $c=d-a+t+1$ and $z_{ij}=z_{[a+l]}$.
In particular,~\eqref{inclutions for Omega c +} holds for $c=d-b+1$ and $z_{ij}=z_{[b]}$.
So, by the second equation in~\eqref{a and b behavior} we have that~\eqref{inclutions of Omega +} holds for
$z_{ij}=z_{[b]}$.
\item $z_{[b]}<0, z_{[b-1]}<0, \ldots, z_{[a]}<0$ and for any $t=1,\ldots,b-a$,~\eqref{inclutions for Omega c +} holds
for $c=d-(b-t)$ and $z_{ij}=z_{[b-t]}$.
In particular,~\eqref{inclutions for Omega c +} holds for $c=d-a$ and $z_{ij}=z_{[a]}$.
So, by the f\/irst equation in~\eqref{a and b behavior} we have that~\eqref{inclutions of Omega +} holds for $z_{ij}=z_{[a]}$. \hfill{$\qed$}
\end{enumerate}
\renewcommand{\qed}{}
\end{proof}

\begin{Definition}
Given $T(Q)$ and $T(R)$ in $\mathcal{B}(T(L))$, we write $T(R)\preceq_{(1)}T(Q)$ if there exist $g\in\mathfrak{gl}(n)$
such that $T(Q)$ appears with nonzero coef\/f\/icient in the decomposition of $g\cdot T(R)$ into a~linear combination of
tableaux.
For any $p\geq 1$ we write $T(R)\preceq_{(p)} T(Q)$ if there exist tableaux $T(L^{(1)})$,\dots, $T(L^{(p)})$, such that
\begin{gather*}
T(R)=T\big(L^{(0)}\big)\preceq_{(1)}T\big(L^{(1)}\big)\preceq_{(1)}\cdots\preceq_{(1)}T\big(L^{(p)}\big)=T(Q).
\end{gather*}
\end{Definition}

As an immediate consequence of the def\/inition of $\preceq_{(p)}$ we have the following.
\begin{Lemma}
\label{properties of relation less or equal}
If $T(Q)$, $T(Q^{(0)})$, $T(Q^{(1)})$ and $T(Q^{(2)})$ are tableaux in $\mathcal{B}(T(L))$ then:
\begin{itemize}\itemsep=0pt
\item[$(i)$] $T(Q^{(0)})\preceq_{(p)} T(Q^{(1)})$ and $T(Q^{(1)})\preceq_{(q)} T(Q^{(2)})$ imply
$T(Q^{(0)})\preceq_{(p+q)} T(Q^{(2)})$;
\item[$(ii)$] $T(Q)\preceq_{(1)} T(Q)$.
\end{itemize}
\end{Lemma}

\begin{Corollary}\label{Gelfand--Tsetlin modules property}\sloppy
If $T(R), T(Q)\in \mathcal{B}(T(L))$ are generic Gelfand--Tsetlin tableaux such that
\mbox{$T(R)\preceq_{(p)} T(Q)$} for some $p\in\mathbb{Z}_{\geq 0}$, then $T(Q)\in U\cdot T(R)$.
\end{Corollary}
\begin{proof}
By Lemma~\ref{elements of Gamma separates tableaux} and the def\/inition of the relation $\preceq_{(1)}$, we f\/irst verify
that ${T(R)\preceq_{(1)}\! T(Q)}$ implies $T(Q)\in U\cdot T(R)$.
Now, using Lemma~\ref{properties of relation less or equal}(i), if $T(R)\preceq_{(p)} T(Q)$ for some~$p$ then $T(Q)\in
U\cdot T(R)$.
\end{proof}

The next theorem provides a~convenient basis for the submodule of $V(T(L))$ generated by a~f\/ixed tableau.
Recall the def\/inition of $\mathcal{N}(T(R))$ in Notation~\ref{notation}(iii)(d).

\begin{Theorem}
\label{basis of submodule generated by T}
For any tableau $T(R)\in\mathcal{B}(T(L))$, $U\cdot T(R)= W(T(R))$.
In particular, $\mathcal{N}(T(R))$ forms a~basis of $U\cdot T(R)$, and the action of $\mathfrak{gl}(n)$ on $U\cdot T(R)$
is given by the Gelfand--Tsetlin formulas.
\end{Theorem}

\begin{proof}
By Proposition~\ref{submodules of V(T(L))}, $U\cdot T(R)\subseteq W(T(R))$.
To prove that $ W(T(R))\subseteq U\cdot T(R)$ we will show that $T(Q)\in U\cdot T(R)$ for any
$T(Q)\in\mathcal{N}(T(R))$.
By Corollary~\ref{Gelfand--Tsetlin modules property}, it is enough to prove that $T(R)\preceq_{(p)} T(Q)$ for some
positive integer~$p$.

Suppose that $T(Q)=T(R+z)\in \mathcal{N}(T(R))$ for some $z\in\mathbb{\mathbb{Z}}^{\frac{n(n-1)}{2}}$.
Let~$t$ be the number of non-zero components of~$z$.
We will prove that $T(R)\preceq_{(p)} T(Q)$ using induction on~$t$.

Let us f\/irst consider the case $t=1$ (the case $t=0$ is trivial, since then $T(Q)=T(R)$) and $z_{ij}>0$.
We will f\/irst prove that $T(R+l\delta^{ij})\preceq_{(1)}T(R+(l+1)\delta^{ij})$ for any $0\leq l \leq z_{ij}-1$.
This will imply
\begin{gather*}
T(R)\preceq_{(1)}T\big(R+\delta^{ij}\big)\preceq_{(1)}T\big(R+2\delta^{ij}\big)\preceq_{(1)}\cdots\preceq_{(1)}T\big(R+z_{ij}\delta^{ij}\big)=T(Q),
\end{gather*}
and then $T(R)\preceq_{(z_{ij})}T(Q)$.
To prove that $T(R+l\delta^{ij})\preceq_{(1)}T(R+(l+1)\delta^{ij})$ we show that the coef\/f\/icient of
$T(R+(l+1)\delta^{ij})$ in the decomposition of $E_{i,i+1}(T(R+l\delta^{ij}))$ is not zero.
In fact, by the Gelfand--Tsetlin formulas, that coef\/f\/icient is
\begin{gather*}
a_{l}:=-\frac{\prod\limits_{k=1}^{i+1}(r_{ij}-r_{i+1,k}+l)}{\prod\limits_{k\neq j}^{i}(r_{ij}-r_{ik}+l)}.
\end{gather*}
Assume that $a_{l}=0$.
Then $r_{ij}-r_{i+1,k}+l=0$ for some~$k$, which implies $\omega_{i+1,k,j}(T(R))$ $=r_{i+1,k}-r_{ij}=l\in\mathbb{Z}_{\geq
0}$.
But, since $T(Q)\in\mathcal{N}(T(R))$, we have{\samepage
\begin{gather*}
l-z_{ij}=r_{i+1,k}-r_{ij}-z_{ij}=\omega_{i+1,k,j}(T(Q))\in\mathbb{Z}_{\geq{0}}.
\end{gather*}
Therefore we have 
$ 0\leq l \leq z_{ij}-1$ and $z_{ij}\leq l$, which is a~contradiction.
Hence,
$T(R)\preceq_{(z_{ij})}T(Q)$.}

Let now $t=1$ and $z_{ij}<0$.
Using the same arguments as in the case $z_{ij}>0$, we prove that $T(R)\preceq_{(-z_{ij})}T(Q)$ using $|z_{ij}|$
applications of $E_{i+1,i}$.
This completes the proof for $t=1$.

Assume now that for any $w\in\mathbb{Z}^{\frac{n(n-1)}{2}}$ with at most~$t$ nonzero components, and such that
$\Omega^{+}(T(R))\subseteq\Omega^{+}(T(R+w))$, we have $T(R)\preceq_{(p)} T(R+w)$ for some~$p$.
Let us consider~$z$ with $t+1$ nonzero components.
Since $\Omega^{+}(T(R))\subseteq \Omega^{+}(T(R+z))$, by Lemma~\ref{paths for connected tableaux}, there exist $i$, $j$
such that
\begin{gather*}
\Omega^{+}(T(R))\subseteq\Omega^{+}\big(T\big(R+z_{ij}\delta^{ij}\big)\big)\subseteq \Omega^{+}(T(R+z)).
\end{gather*}
Using the induction hypothesis for the pairs of tableaux $(T(R), T(R+z_{ij}\delta^{ij}))$ and
$(T(R+z_{ij}\delta^{ij}),T(R+z))$, there exist $p,q\in\mathbb{Z}_{\geq 0}$ such that
\begin{gather*}
T(R)\preceq_{(p)} T\big(R+z_{ij}\delta^{ij}\big)\qquad \text{and} \qquad T\big(R+z_{ij}\delta^{ij}\big)\preceq_{(q)} T(R+z).
\end{gather*}
Thus, by Lemma~\ref{properties of relation less or equal}(i), $T(R)\preceq_{(p+q)}T(R+z)$.
\end{proof}

\begin{Proposition}
\label{when two tableaux generate the same submodule}
Let $T(R)$ and $T(Q)$ be in $\mathcal{B}(T(L))$.
Then $U\cdot T(R)= U\cdot T(Q)$ if and only if $\Omega^{+}(T(Q))=\Omega^{+}(T(R))$.
\end{Proposition}
\begin{proof}\sloppy
Using Theorem~\ref{basis of submodule generated by T} and the def\/initions of $ W(T(R))$, $W(T(Q))$, $\Omega^{+}(T(R))$,
and $\Omega^{+}(T(Q))$, we can prove a~stronger statement: $U\cdot T(R) \subseteq U\cdot T(Q)$ if and only if
$\Omega^{+}(T(Q))\subseteq\Omega^{+}(T(R))$.
\end{proof}

\begin{Corollary}
$U \cdot T(R)=V(T(L))$ whenever $\Omega^{+}(T(R))=\varnothing$.
\end{Corollary}

\begin{Definition}
We will write $T(Q)\sim_{\Omega^+}T(R)$ if $\Omega^{+}(T(R)) = \Omega^{+}(T(Q))$.
\end{Definition}

\begin{Proposition}
Every submodule of $V(T(L))$ is finitely generated.
\end{Proposition}
\begin{proof}
Let~$N$ be any submodule of $V(T(L))$ and~$\Phi$ the set of all tableaux $T(R)$ in~$N$ such that
$\Omega^{+}(T(P))\subseteq\Omega^{+}(T(R))$ implies $\Omega^{+}(T(P))=\Omega^{+}(T(R))$.
By Theorem~\ref{basis of submodule generated by T}, $N=\sum\limits_{T(R)\in\Phi} U\cdot T(R)$ and by
Proposition~\ref{when two tableaux generate the same submodule}, we can write $N=\bigoplus_{T(R)\in\tilde{\Phi}} U\cdot
T(R)$, where $\tilde{\Phi}$ is a~set of distinct representatives of $\Phi/\sim_{\Omega^+}$ (hence
$\Omega^{+}(T(R))\neq\Omega^{+}(T(Q))$ for any $T(R)$, $T(Q)$ in $\tilde{\Phi}$).
Now, since $\Omega(T(L))$ is a~f\/inite set, then $\tilde{\Phi}$ is f\/inite.
\end{proof}

\subsection{Basis for irreducible modules containing a~given tableau}
By Theorem~\ref{basis of submodule generated by T}, the module generated by a~tableau $T(R)$ has basis
$\mathcal{N}(T(R))$.
For the purpose of the next theorem let us introduce the following equivalence on ${\mathbb C}^{\frac{n(n+1)}{2}}$.

\begin{Definition}
We write $z\sim w$ for $z,w \in {\mathbb C}^{\frac{n(n+1)}{2}}$ if and only if one of the two cases hold.
\begin{itemize}\itemsep=0pt
\item[(i)] $z-w \in {\mathbb Z}^{\frac{n(n-1)}{2}}$ and $z \sim_{\Omega^+} w$.
\item[(ii)] $z \in G w$.
\end{itemize}
\end{Definition}

Now we are ready to formulate and prove the main theorem in the paper.

\begin{Theorem}
\label{characterization of irreducible basis}
The irreducible module containing $T(R)$ has a~basis of tableaux
\begin{gather*}
\mathcal{I}(T(R))=\big\{T(Q)\in\mathcal{B}(T(R)): \Omega^{+}(T(Q))=\Omega^{+}(T(R))\big\}.
\end{gather*}
The action of $\mathfrak{gl}(n)$ on this irreducible module is given by the Gelfand--Tsetlin
formulas~\eqref{Gelfand--Tsetlin formulas}.
Therefore the set of irreducible generic Gelfand--Tsetlin modules is in one-to-one correspondence with $ {\mathbb
C}^{\frac{n(n+1)}{2}}_{\rm gen}/ \sim$, where ${\mathbb C}^{\frac{n(n+1)}{2}}_{\rm gen}$ stands for the set of generic
vectors in ${\mathbb C}^{\frac{n(n+1)}{2}}$.
\end{Theorem}

\begin{proof}

For each tableau $T(R)$, we have an explicit construction of the module containing $T(R)$ (recall
Def\/inition~\ref{def-cont}):
\begin{gather*}
M(T(R)):=U\cdot T(R)/\left(\sum U\cdot T(Q)\right),
\end{gather*}
where the sum is taken over tableaux $T(Q)$ such that $T(Q)\in U\cdot T(R)$ and $U\cdot T(Q)$ is a~proper submodule of
$U\cdot T(R)$.

The module $M(T(R))$ is simple.
Indeed, this follows from the fact that for any nonzero tableau $T(S)$ in $M(T(R))$ we have $U\cdot T(S)=U\cdot T(R)$
and, hence, $T(S)$ generates $M(T(R))$.

By Theorem~\ref{basis of submodule generated by T} and Proposition~\ref{when two tableaux generate the same submodule},
a~basis for a~proper submodule $U\cdot T(Q)$ of $U\cdot T(R)$ is
$\{T(S):\Omega^{+}(T(R))\subsetneq\Omega^{+}(T(Q))\subseteq \Omega^{+}(T(S))\}$ so, a~basis for the module $\sum U\cdot
T(Q)$ is $\{T(S):\Omega^{+}(T(R))\subsetneq\Omega^{+}(T(S))\}$.
Therefore, $\mathcal{I}(T(R))$ is a~basis for $M(T(R))$.

To show that ${\mathbb C}^{\frac{n(n+1)}{2}}_{\rm gen}/ \sim$ parameterizes the set of all irreducible generic
Gelfand--Tsetlin modules we use Theorem~\ref{uniqueness in generic case} and the fact that $\ell, \ell'\in\Specm\Lambda$ lie over the same $\mathsf m$
in $\Specm\Gamma$ if and only if $\ell\in G\ell'$ (see Remark~\ref{maximal ideals differ by permutations}).
\end{proof}

\section{Number of irreducible modules in generic blocks}

\begin{Definition}
For any generic tableau $T(L)$, the {\em block associated with $T(L)$} is the set of all Gelfand--Tsetlin
$\mathfrak{gl}(n)$-modules with Gelfand--Tsetlin support contained in $\operatorname{Supp}_{\rm GT}(V(T(L)))$.
\end{Definition}

Theorem~\ref{characterization of irreducible basis} describe explicit bases of the irreducible modules in the block
associated with $V(T(L))$.
In this section we will use this description to compute the number of nonisomorphic irreducible modules in this block.

\begin{Definition}
For any $T(R)=T(r_{ij})\in\mathcal{B}(T(L))$, $1<p\leq n$ and $1\leq u\leq p-1$, def\/ine $d_{pu}(T(R))$ to be the number
of distinct elements in
\begin{gather*}
\{r_{ps}:(p,s,u)\in\Omega(T(R))\}.
\end{gather*}
\end{Definition}

\begin{Remark}
For any generic tableau $T(R)=T(r_{ij})\in\mathcal{B}(T(L))$ of height~$n$ we have:
\begin{enumerate}[(i)]\itemsep=0pt
\item $d_{pu}(T(L))=d_{pu}(T(R))$ for any $1<p\leq n$, $1\leq u\leq p-1$;
\item if $p\neq n$, then $d_{pu}(T(R))\leq 1$ for any $1\leq u\leq p-1$.
\end{enumerate}
\end{Remark}

\begin{Example}
Suppose $a,b,c\in\mathbb{C}$ are such that $\{a-b,a-c,b-c\}\cap \mathbb{Z}=\varnothing$.
If $R=(a,a-1,b|a,b|c)$, then
\begin{center}
\hspace{1.5cm}\Stone{$a$}\Stone{$a-1$}\Stone{$b$}
\\[0.2pt] $T(R)$:=\hspace{0.2cm} \Stone{$a$}\Stone{$b$}
\\[0.2pt] \hspace{1.4cm}\Stone{$c$}
\\
\end{center}
$d_{31}(T(R))=2$, $d_{32}(T(R))=1$, $d_{21}(T(R))=0$ and $d_{22}(T(R))=0$.
\end{Example}

\begin{Remark}
\label{one to one correspondence}
For each tableau $T(R)$ we have an one-to-one correspondence between the set $\{0,1,\ldots,d_{pu}(T(L))\}$ and the
subset $\{0,i_{1},\ldots,i_{d_{pu}(T(L))}\}$ of $\{0,1,\ldots,p\}$ def\/ined as follows: \mbox{$i_{1}=1$} and $i_{k} = \min \{x:
r_{px}\notin\{r_{pi_{1}},\ldots,r_{pi_{k-1}} \} \}$.
\end{Remark}

\begin{Theorem}
For any generic tableau $T(L)$, the number of irreducible modules in the block associated with $T(L)$ is
\begin{gather*}
\prod\limits_{1\leq u\leq p-1< n}(d_{pu}(T(L))+1).
\end{gather*}
In particular, $V(T(L))$ is irreducible if and only if $d_{pu}(T(L))=0$ for any~$p$ and~$u$, or equivalently, if and
only if $\Omega(T(L))=\varnothing$.
\end{Theorem}
\begin{proof}
By Theorem~\ref{characterization of irreducible basis}, the irreducible modules are in one-to-one correspondence with
the subsets of $\Omega(T(L))$ of the form $\Omega^{+}(T(L+z))$.
For any $T(R)\in\mathcal{B}(T(L))$, we can decompose $\Omega(T(R))$ into a~disjoint union
$\Omega(T(R))=\bigsqcup_{p,u}\Omega_{pu}(T(R))$, where
\begin{gather*}
\Omega_{p,u}(T(R))=\{(p,1,u),(p,2,u),\ldots,(p,p,u)\}\cap \Omega(T(R)).
\end{gather*}
Now, if $\Omega^{+}_{p,u}(T(R)):=\Omega_{p,u}\cap\Omega^{+}(T(R))$, one can write
$\Omega^{+}(T(R))=\bigsqcup_{p,u}\Omega^{+}_{pu}(T(R))$.
For $p$, $u$ f\/ixed, let us denote by $s_{p,u}$ the number of dif\/ferent subsets of the form $\Omega^{+}_{p,u}(T(R))$.
So, the number of dif\/ferent subsets of the form $\Omega^{+}(T(R))$ is $\prod\limits_{p,u}s_{p,u}$.

Let $\{T(R^{(i)})\}_{i=1}^{s_{pu}}$ be a~set of tableaux such that $\{\Omega^{+}_{p,u}(T(R^{(i)}))\}_{i=1}^{s_{pu}}$ is
the set of all distinct sets of the form $\Omega^{+}_{p,u}(T(R))$.
We have a~one-to-one correspondence between $\{T(R^{(i)})\}_{i=1}^{s_{pu}}$ and the set
$\{0,i_{1},\ldots,i_{d_{pu}(T(L))}\}$ constructed as in Remark~\ref{one to one correspondence}.
More explicitly, this correspondence is def\/ined my the map:
\begin{gather*}
T(R^{(i)})\to
\begin{cases}
\min\{j: (p,j,u)\in\Omega^{+}(T(R^{(i)}))\}, & \text{if} \ \ \Omega^{+}_{pu}(T(R^{(i)}))\neq \varnothing,
\\
0, & \text{if} \ \ \Omega^{+}_{pu}(T(R^{(i)}))= \varnothing.
\end{cases}
\end{gather*}
Therefore, $s_{pu}=d_{pu}(T(L))+1$.
\end{proof}

\subsection*{Acknowledgements}

V.F. gratefully acknowledges the hospitality and excellent working conditions
at the CRM, University of Montreal, where part of this work was completed.
V.F.~is supported in part by the CNPq grant (301320/2013-6) and by the Fapesp grant (2014/09310-5).
D.G.~is supported in part by the Fapesp grant (2011/21621-8) and by the NSA grant H98230-13-1-0245.
L.E.R.~is supported by the Fapesp grant (2012/23450-9).

\pdfbookmark[1]{References}{ref}
\LastPageEnding

\end{document}